\documentclass[11pt]{amsart}
\theoremstyle{plain}
\newtheorem{thm}{Theorem}[section]
\newtheorem{theorem}[thm]{Theorem}

\theoremstyle{definition}

\newtheorem{definition}[thm]{Definition}

\numberwithin{equation}{section}

\newcommand{\T}{{\mathbb T}}

\title[Pseudo-isomorphisms and Monge-Ampere equations]{Pseudo-isomorphisms in dimension $3$ and applications to complex Monge-Ampere equation}

\author{Tuyen Trung Truong}

\address{Department of Mathematics, Syracuse University, Syracuse NY 13244} \email{tutruong@syr.edu}

\thanks{}

\begin{document}

\maketitle

\begin{abstract}
Let $X$ and $Y$ be compact K\"ahler manifolds of dimension $3$. A bimeromorphic map $f:X\rightarrow Y$ is pseudo-isomorphic if $f:X-I(f)\rightarrow Y-I(f^{-1})$ is an isomorphism.  

In this paper we investigate some properties of pseudo-isomorphisms. As an application, we associate to any pseudo-isomorphism in dimension $3$ and a smooth closed $(3,3)$ form $\delta$ on $X\times X$ representing the cohomology class of the diagonal $\Delta _X$, a Monge-Ampere operator $MA(f^*(\theta ),\delta )=f^*(\theta )\wedge f^*(\theta )\wedge f^*(\theta )$, here $\theta$ is a smooth closed $(1,1)$ form on $Y$. We show that this Monge-Ampere operator  is independent of the choice of $\delta$, if the following cohomologous condition is satisfied:

{\bf Condition.} For any curve $C\subset I(f^{-1})$, we have $\{\theta \}.\{C\}=0$ in cohomology. 

We conclude the paper examining a simple pseudo-isomorphism in dimension $3$.  

\end{abstract}

\section{Introduction}

Let $X$ and $Y$ be compact K\"ahler manifolds of dimension $3$. A bimeromorphic map $f:X\rightarrow Y$ is pseudo-isomorphic if the map $g=f|_{X-I(f)}:X-I(f)\rightarrow Y-I(f^{-1})$ is an isomorphism. We let $\Gamma _g\subset (X-I(f))\times (Y-I(f^{-1}))$ be the graph of $g$, and $\Gamma _f=$ the closure of $\Gamma _g$ in $X\times Y$ the graph of $f$. Let $\pi _1,\pi _2:~\Gamma _f,\Gamma _g\rightarrow X,Y$ be the natural projections. 

For a meromorphic map, Meo \cite{meo} defined the pullback of a positive closed $(1,1)$ current and a quasi-psh function. In \cite{truong1}, we showed that for a pseudo-isomorphism in dimension $3$, we can pullback and pushforward any positive closed $(2,2)$ currents. Moreover, these pullback  and pushforward are continuous with the weak topology on currents.    

In this paper, we investigate further properties of pseudo-isomorphisms and give applications to complex Monge-Ampere operators in dimension $3$.

 Our first results are the following. 

\begin{theorem}
Let $f:X\rightarrow Y$ be a pseudo-isomorphism in dimension $3$. 

Let $T_1$ and $T_2$ be differences of positive closed $(1,1)$ currents which are smooth outside a curve. Then the currents $f^*(T_1)\wedge f^*(T_2)$ and $f^*(T_1\wedge T_2)$ are well-defined. (The first one being by dimension reason, the second one following \cite{truong0, truong1}. Please see the proof of the theorem for more details.)

The difference current $f^*(T_1)\wedge f^*(T_2)-f^*(T_1\wedge T_2)$, which has support in $I(f)$, depends only on the cohomology classes of $T_1$ and $T_2$.

Moreover, if for every curve $C\subset I(f^{-1})$ we have in cohomology: $\{T_1\}.\{C\}=0$, then $f^*(T_1)\wedge f^*(T_2)-f^*(T_1\wedge T_2)=0$.
\label{Theorem0}\end{theorem}

\begin{theorem}
Let $u$ be a quasi-psh function on $X$ and $T$ a positive closed $(1,1)$ current on $Y$ which is smooth on $Y-C$ where $C$ is a curve which is not contained in $I(f^{-1})$. Then the current $uf^*(T)$ is well-defined on $X$. Moreover, the current $\pi _1^*(uf^*(T))$ has bounded mass on $\Gamma _g$.
\label{Theorem1}\end{theorem}

By Theorem \ref{Theorem1}, we can extend the current $\pi _1^*(uf^*(T))$ on $\Gamma _g$ by zero to $\Gamma _f$. We denote this current by $\widetilde{(\pi _1^*(uf^*(T ))\wedge [\Gamma _g])}$. We then define for any smooth function $\phi$ on $X$
\begin{equation}
f_*(\varphi dd^cuf^*(T)):=(\pi _2)_*(\pi _1^*(\varphi )dd^c\widetilde{(\pi _1^*(uf^*(T))\wedge [\Gamma _g])}).
\label{Equation1}\end{equation}

In case $u$ is a smooth function, the definition in Equation (\ref{Equation1}) may be different from the correct value
\begin{equation}
(\pi _2)_*(\pi _1^*(\varphi )dd^c(\pi _1^*(uf^*(T))\wedge [\Gamma _f])).
\label{Equation2}\end{equation}

However, under one cohomologous condition, Equations (\ref{Equation1}) and (\ref{Equation2}) are the same.   
\begin{theorem}
Assumptions are as in Theorem \ref{Theorem1}. Moreover, assume that $u$ is smooth, and either 

(i) $dd^cu=0$, 

or

(ii) for any curve $C$ in the indeterminate set $I(f^{-1})$, we have in cohomology $\{T\}.\{C\}=0$. 

Then Equations \ref{Equation1} and \ref{Equation2} are the same. 
\label{Theorem2}\end{theorem}

As an application, we define the Monge-Ampere operators $MA(f^*(\theta ))=f^*(\theta )\wedge f^*(\theta )\wedge f^*(\theta )$, where $\theta$ is a smooth closed $(1,1)$ form, as follows. We write $f^*(
\theta )=\Omega +dd^cu$, where $\Omega $ is a smooth closed $(1,1)$ form and $u$ is a quasi-psh function on $X$. If $\varphi$ is a smooth function on $X$ then we define 
\begin{equation}
<f^*(\theta )\wedge f^*(\theta )\wedge f^*(\theta ),\varphi >:=\int _Y\varphi f^*(\theta )\wedge f^*(\theta )\wedge \Omega +\int _Y\theta \wedge f_*(\varphi dd^cuf^*(\theta )).
\label{Equation3}\end{equation}
Here $f_*(\varphi dd^cuf^*(\theta ))$ is defined in Equation (\ref{Equation1}).

{\bf Remark 1.} The Monge-Ampere defined above is correct in cohomology, that is the total mass of $MA(f^*(\theta ))$ is the same as the intersection in cohomology $\{f^*(\theta )\}.\{f^*(\theta )\}.\{f^*(\theta )\}$.

\begin{theorem}
For a fixed choice of $\Omega$, the Monge-Ampere operator defined in Equation (\ref{Equation3}) is independent of the choice of the quasi-potential $u$.

Moreover, the Monge-Ampere operator is independent of $\Omega$, if the following cohomologous condition is satisfied: for every curve $C\subset I(f^{-1})$, then in cohomology $\{\theta \}.\{C\}=0$.
\label{Theorem3}\end{theorem}

{\bf Remark 2.} Given a smooth closed $(3,3)$ form $\delta$ on $X\times X$ which is cohomologous to the diagonal $\Delta _X$, we have a uniform choice of a smooth representing $\Omega _{\delta ,T}$ of positive closed $(1,1)$ currents $T$ on $X$ by the following formula
\begin{eqnarray*}
\Omega _{\delta ,T}:=(\pi _1)_*(\pi _2^*(T)\wedge \delta ).
\end{eqnarray*}
Here $\pi _1,\pi _2:X\times X\rightarrow X$ are the natural projections.

Therefore, we have the following definition
\begin{definition}
Let $f:X\rightarrow Y$ be a pseudo-isomorphism in dimension $3$. Let $\delta$ be a smooth closed $(3,3)$ form on $X$ having the same cohomology class as that of the diagonal $\Delta _X$. 

We define the Monge-Ampere operator $MA(f^*(\theta ),\delta )$ by the following formula
\begin{eqnarray*}
<MA(f^*(\theta ),\delta ),\varphi >:=\int _Y\varphi f^*(\theta )\wedge f^*(\theta )\wedge \Omega _{\delta ,f^*(\theta )} +\int _Y\theta \wedge f_*(\varphi dd^cuf^*(\theta ))
\end{eqnarray*}
Here $f^*(\theta )=\Omega _{\delta ,f^*(\theta )}+dd^cu$. 

By Theorem \ref{Theorem3}, this is independent of the choice of $u$.
\label{Definition1}\end{definition}

In the last Section, we will apply Theorem \ref{Theorem3} to a specific pseudo-automorphism $J_X$ in dimension $3$. We obtain the following result
\begin{theorem}
Let $J_X:X\rightarrow X$ be the pseudo-automorphism considered in Section \ref{Section1}. 

1) There is a unique non-zero cohomology class $\eta \in H^{1,1}(X)$ such that if $\theta$ is a smooth closed $(1,1)$ form for which the Monge-Ampere operator $MA(.,\delta )$ in Definition \ref{Definition1} is independent of the choice of $\delta$, then in cohomology $\{\theta\}$ is a multiple of $\{\eta\}$.

Moreover, the class $\eta$ is nef.

2) There is a positive closed smooth $(1,1)$ form $\theta$ on $X$ such that for every choice of $\delta$ then the Monge-Ampere operator $MA(f^*(\theta ),\delta )$ is a signed measure with support in $I(f)$ and with total mass $-3$.
\label{Theorem4}\end{theorem}

{\bf Acknowledgments.} The author would like to thank Eric Bedford, Tien-Cuong Dinh, Mattias Jonsson, Nessim Sibony and Yuan Yuan for their crucial comments. The author also would like to thank Muhammed Alan for helpful discussions on the topic.

\section{Proofs of Theorems \ref{Theorem0}, \ref{Theorem1} and \ref{Theorem2}}
\label{Section0}

\subsection{Proof of Theorem \ref{Theorem0}}

Let $Z$ be a resolution of singularity of the graph $\Gamma _f$ of $f$, and let $\pi ,h:Z\rightarrow X,Y$ be the induced holomorphic maps to the source and the target. We can choose $\pi $ such that it is a composition of blowups at a smooth point or curve $\pi =\pi _m\circ \circ \ldots \pi _2 \circ \pi _1$, and moreover, the images by $\pi$ and $h$ of the exceptional divisor of $\pi$ are contained in $I(f)$ and $I(f^{-1})$. 

1) We first prove part 1) in the case $T_1$ and $T_2$ are smooth closed $(1,1)$ forms. Then 
\begin{eqnarray*}
f^*(T_1)=\pi _*h^*(T_1),~f^*(T_2)=\pi _*h^*(T_2),~f^*(T_1\wedge T_2)=\pi _*(h^*(T_1)\wedge h^*(T_2)). 
\end{eqnarray*}

Here $\pi _*=(\pi _1)_*\circ \ldots (\pi _m)_*$. Define Since $\alpha =h^*(T_1)$ and $\beta =h^*(T_2)$. These are  smooth closed $(1,1)$ forms on $Z$.

a) First, consider the blowup $\pi _1$. If $\pi _1$ is blowup at a point then there is nothing to prove. Hence we consider the case $\pi _1$ is the blowup at a smooth curve $D_1$. Let $F_1$ be a fiber of the restriction of $\pi _1$ to the exceptional divisor $\pi _1^{-1}(D_1)$.  by Lemma 4 in \cite{truong2} we have
\begin{equation}
(\pi _1)_*((\pi _1)^*(\pi _1)_*(\alpha )\wedge \beta )-(\pi _1)_*(\alpha \wedge \beta )=\{\alpha .F_1\}\{\beta .F_1\}[D_1].
\label{Equation4}\end{equation}  
Therefore the current $(\pi _1)_*((\pi _1)^*(\pi _1)_*(\alpha )\wedge \beta )-(\pi _1)_*(\alpha \wedge \beta )$ depends only on the cohomology classes of $\alpha $ and $\beta$, therefore depends only on the cohomology classes of $T_1$ and $T_2$. 

We now claim that the currents $(\pi _1)_*((\pi _1)^*(\pi _1)_*(\alpha )\wedge \beta )$ and $(\pi _1)_*(\alpha )\wedge (\pi _1)_*(\beta )$ are the same. In fact, the currents $(\pi _1)_*(\alpha )$ and $(\pi _1)_*(\beta )$ are   are differences of positive closed $(1,1)$ currents which are smooth outside a curve, so the intersection $(\pi _1)_*(\alpha )\wedge (\pi _1)_*(\beta )$ is well-defined (see Section 4, Chapter 3 in \cite{demailly}). If we consider a smooth approximation $S_n$ of $(\pi _1)_*(\alpha )$, then $(\pi _1)^*(S_n)\wedge \beta $ converges to $(\pi _1)^*(\pi _1)_*(\alpha )\wedge \beta$, therefore
\begin{eqnarray*}
(\pi _1)_*((\pi _1)^*(\pi _1)_*(\alpha )\wedge \beta )=\lim _{n\rightarrow\infty}(\pi _1)_*(\pi _1^*(S_n)\wedge \beta )=\lim _{n\rightarrow\infty}S_n\wedge (\pi _1)_*(\beta )=(\pi _1)_*(\alpha )\wedge (\pi _1)_*(\beta ).
\end{eqnarray*}
The second equality follows from the projection formula, and the third equality follows if we choose $S_n$ good enough. 

Therefore, we showed that the current
\begin{eqnarray*}
(\pi _1)_*(\alpha )\wedge (\pi _1)_*(\beta )-(\pi _1)_*(\alpha \wedge \beta )
\end{eqnarray*}
depends only on the cohomology classes of $T_1$ and $T_2$.

b) We consider the second blowup $\pi _2$. Here we can not apply directly Lemma 4 in \cite{truong2} since now $(\pi _1)_*(\alpha )$ and $(\pi _1)_*(\beta )$ are not smooth. However, we can approximate $(\pi _1)_*(\alpha )$ by smooth closed $(1,1)$ forms $S_n$ and $(\pi _1)_*(\beta )$ by smooth closed $(1,1)$ forms $R_n$. Then each term 
\begin{eqnarray*}
(\pi _2)_*(S_n)\wedge (\pi _2)_*(R_n)-(\pi _2)_*(S_n\wedge R_n)
\end{eqnarray*}
depends only on the cohomology classes of $S_n$ and $R_n$. Hence the limit depends only on the cohomology classes of $T_1$ and $T_2$. We now show that $(\pi _2)_*(S_n)\wedge (\pi _2)_*(R_n)$ converges to $(\pi _2)_*(\pi _1)_*(\alpha ) \wedge (\pi _2)_*(\pi _1)_*(\beta )$. Note that again, $(\pi _2)_*(\pi _1)_*(\alpha )$ and $(\pi _2)_*(\pi _1)_*(\beta )$ are smooth outside a curve, so their wedge product is well-defined. We can choose $S_n$ with the following properties: $S_n=\Omega _n+dd^cu_n$, where $\Omega $ are smooth closed $(1,1)$ form converging uniformly to a $C^2$ closed (1,1) form, and $u_n$ are differences of quasi-psh functions converging locally uniformly outside a curve. We can choose $R_n$ such that it converges locally uniformly outside a curve. Then an argument similar to that in proof of Lemma 5 in \cite{truong1} shows that the sequence $u_nR_n$ converges (the sequence converges locally uniformly outside a curve, and the limit has no mass on curves because of the dimension reason). It follows that $S_n\wedge R_n=\Omega _n\wedge R_n+dd^c(u_nR_n)$ converges, and the limit is exactly $(\pi _2)_*(\pi _1)_*(\alpha ) \wedge (\pi _2)_*(\pi _1)_*(\beta )$. 

By a similar argument, we can show that $(\pi _2)_*(S_n\wedge R_n)$ converges to $(\pi _2)_*((\pi _1)_*(\alpha )\wedge (\pi _1)_*(\beta ))$. Combining with Step a) we conclude that the current
\begin{eqnarray*}
(\pi _2)_*(\pi _1)_*(\alpha ) \wedge (\pi _2)_*(\pi _1)_*(\beta )-(\pi _2)_*(\pi _1)_*(\alpha \wedge \beta )
\end{eqnarray*}
depends only on the cohomology classes of $T_1$ and $T_2$.

c) Working inductively on the number of blowups, we conclude that for smooth closed $(1,1)$ forms $T_1$ and $T_2$
\begin{eqnarray*}
f^*(T_1)\wedge f^*(T_2)-f^*(T_1\wedge T_2)
\end{eqnarray*}
depends only on the cohomology classes of $T_1$ and $T_2$. 

2) Consider the general case. Since $T_1$ and $T_2$ are smooth outside a curve, so is $f^*(T_1)$ and $f^*(T_2)$. Therefore $f^*(T_1)\wedge f^*(T_2)$ is well-defined. Since $T_1\wedge T_2$ is a positive closed $(2,2)$ current, $f^*(T_1\wedge T_2)$ is well-defined as shown in \cite{truong0, truong1}. Moreover, if we approximate $T_1$ and $T_2$ by appropriate smooth closed $(1,1)$ forms $S_n$ and $R_n$ (for examples like the ones used in the proof of Step b) ), then $f^*(S_n\wedge R_n)$ converges to $f^*(T_1\wedge T_2)$ and $f^*(S_n)\wedge f^*(R_n)$ converges to $f^*(S_n)\wedge f^*(R_n)$. Here we use that since $f$ is a pseudo-automorphism in dimension $3$, if $Sn$ converges locally uniformly outside a curve, then so is $f^*(S_n)$. Also, if $\Omega _n$ converges uniformly to $\Omega $, then $f^*(\Omega _n-\Omega )\wedge f^*(R_n)$ converges to $0$. 

3) We now prove the last assertion of Theorem \ref{Theorem0}. From Equation (\ref{Equation4}), it follows that if in cohomology $\{h^*(\T_1 )\}.\{D\}=0$ for every curve $D$ in the exceptional divisors of the holomorphic map $\pi :Z\rightarrow X$, then $f^*(T_1)\wedge f^*(T_2)-f^*(T_1\wedge T_2)=0$. Since
\begin{eqnarray*}
\{h^*(T_1 )\}.\{D\}=\{T_1\}.h_*\{D\},
\end{eqnarray*}
and $h(D)$ is contained in $I(f^{-1})$, if $\{T_1\}.\{C\}=0$ for every curve $C$ in $I(f^{-1})$ then we also have $\{h^*(T_1 )\}.\{D\}=0$ as wanted.   

\subsection{Proof of Theorem \ref{Theorem1}}

 The pullback of a positive closed $(1,1)$ current and of a quasi-psh function is given by Meo \cite{meo}. Since $f^*(T)$ is smooth outside $f^{-1}(C)$ which is a curve (here we use that $f$ is a pseudo-isomorphism), the current $uf^*(T)$ is well-defined on $X-f^{-1}(C)$ (see e.g. Section 4, Chapter 3 in Demailly \cite{demailly}). Since $f^{-1}(C)$ has codimension at least $2$, we can extend $uf^*(T)$ to be a negative $(1,1)$ current on $X$. Moreover, $dd^c(uf^*(T))$ is a difference of two positive closed $(2,2)$ currents, hence it is DSH as defined in Dinh-Sibony \cite{dinh-sibony1, dinh-sibony2}.  

We now show that the mass of $\pi _1^*(uf^*(T))$ is bounded on $\Gamma _g$. Let $\omega _X$ and $\omega _Y$ be K\"ahler forms on $X$ and $Y$.  Since $\pi _1^*(\omega _X)+\pi _2^*(\omega _Y)$ is a K\"ahler form on $X\times Y$, it suffices to show that each of the following integrals 
\begin{eqnarray*}
\int _{\Gamma _g}\pi _1^*(uf^*(T))\wedge \pi _1^*(\omega _X^2),~\int _{\Gamma _g}\pi _1^*(uf^*(T))\wedge \pi _1^*(\omega _X)\wedge \pi _2^*(\omega _Y),~\int _{\Gamma _g}\pi _1^*(uf^*(T))\wedge \pi _2^*(\omega _Y^2)
\end{eqnarray*}
is bounded. 

The first term 
\begin{eqnarray*}
\int _{\Gamma _g}\pi _1^*(uf^*(T))\wedge \pi _1^*(\omega _X^2)=\int _{X-I(f)}uf^*(T)\wedge \omega _X^2
\end{eqnarray*}
is clearly bounded. This follows from the Oka's principle, see Fornaess-Sibony \cite{fornaess-sibony2}. Here $uf^*(T)$ is a DSH (1,1) current on $X-I(f)$, and $I(f)$ is of dimension $1$ which is smaller than the dimension of $uf^*(T)$. Therefore $uf^*(T)$ extends uniquely as a DSH $(1,1)$ current on $Y$. A similar argument was used in the proofs of Theorem 6 in \cite{truong0} and Lemma 5 in \cite{truong1}. 
 
The second term 
\begin{eqnarray*}
\int _{\Gamma _g}\pi _1^*(uf^*(T))\wedge \pi _1^*(\omega _X)\wedge \pi _2^*(\omega _Y)=\int _{Y-I(f^{-1})}f_*(u\omega _X)\wedge T \wedge \omega _Y
\end{eqnarray*}
is also bounded. Here we use that $f_*(f^*(\theta ))=\theta$ on $Y-I(f^{-1})$, since $f$ is an isomorphism from $X-I(f)$ to $Y-I(f^{-1})$. Let $A\subset I(f^{-1})$ be the finite set where $T$ is not smooth. Since $T$ is smooth on $Y-A-C$ by assumption  and $f_*(\omega _X)$ is smooth outside $I(f^{-1})$, the current $f_*(u\omega _X)\wedge T $   is well-defined on $Y-A$ as a DSH $(2,2)$ current. The Oka's principle again implies that we can extend $f_*(u\omega _X)\wedge T $ uniquely as a DSH current on all of $Y$.  

The last term 
\begin{eqnarray*}
\int _{\Gamma _g}\pi _1^*(uf^*(\theta ))\wedge \pi _2^*(\omega _Y^2)=\int _{Y-I(f^{-1})}f_*(u)T \wedge \omega _Y^2
\end{eqnarray*}
is also bounded. 

\subsection{Proof of Theorem \ref{Theorem2}}

Let $\pi , h:Z\rightarrow X,Y$ be the maps given in the beginning of the proof of Theorem \ref{Theorem0}. Since $f$ is a pseudo-isomorphism, we can choose such that if $D$ is a curve in the exceptional divisor of $\pi$ then $\pi (D)$ is contained in $I(f)$ and $h(D)$ is contained in $I(f^{-1})$. 

We can write 
\begin{eqnarray*}
\pi ^*(f^*(T))=h^*(T)+\sum _{j}\lambda _j[V_j],
\end{eqnarray*}
where $V_j$ are irreducible components of the exceptional divisor of $\pi$, and $\lambda _j$ are constants. 

We note that $h^*(T)=h^o(T)+\sum _{j'}\lambda _{j'}[V_{j'}]$. Here $h^o(T)$ is the strict transform of $T$ hence has no mass on proper analytic sets, and $V_j$ are subvarieties of dimension $2$ in the exceptional divisor of $\pi$ such that $h(V_j')$ is a point. Therefore, for a smooth function $u$, the difference between Equations (\ref{Equation1}) and (\ref{Equation2}) is 
\begin{eqnarray*}
\sum _j\lambda _jh_*(\pi ^*\varphi \wedge \pi ^*dd^cu[V_j]).
\end{eqnarray*}
Therefore, the difference is $0$ when either the restriction of $dd^cu=0$ to the exceptional divisor is $0$, or $\lambda _j=0$ for all $j$.

Assumption (i) of the theorem implies $dd^cu=0$. We now show that condition (ii) is equivalent to $\lambda _j=0$ for all $j$. 

It is known that the cohomology classes of the irreducible components of the exceptional divisor of $\pi$ are linearly independent, see Proposition 5.5 in Dinh-Sibony \cite{dinh-sibony2}. Therefore, to show that $\lambda _j=0$ for all $j$, it suffices to show that $\pi ^*(f^*(T))-h^*(T )=0$ in $H^{1,1}(Z)$. By Poincar\'e duality, we only need to show that for any $\eta \in H^{2,2}(Z)$ we have
\begin{eqnarray*}
\{\pi ^*(f^*(T))-h^*(T)\}.\eta =0.
\end{eqnarray*}

Since $\pi :Z\rightarrow X$ is a composition of blowups, we can write $\eta =\pi ^*\xi +\{D\}$, where $\xi \in H^{2,2}(X)$ and $D$ is a linear combination of curves in the exceptional divisor of $\pi$. 

By the projection formula
\begin{eqnarray*}
\{\pi ^*(f^*(T))-h^*(T)\}.\pi ^*(\xi )=\pi _*\{\pi ^*(f^*(T))-h^*(T)\}.\xi =0.
\end{eqnarray*}
The second equality follows since by definition $\pi _*(h^*(T))=f^*(T)$. 

Now let $D$ be a curve in the exceptional divisor of $\pi$. By the choice of the desingularization $Z$, $\pi (D)$ is contained in $I(f)$ and $h(D)$ is contained in $I(f^{-1})$. Therefore, both $f_*\pi _*\{D\}$ and $h_*\{D\}$ are cohomologus to a linear combination of curves in $I(f^{-1})$. We have by assumption (ii) of the theorem
\begin{eqnarray*}
\{\pi ^*(f^*(T))\}.\{D\}&=&\{T\}.f_*\pi _*\{D\}=0,\\
\{h^*(T)\}.\{D\}&=&\{T\}.h_*\{D\}=0.
\end{eqnarray*}
  
Hence the proof of Theorem \ref{Theorem2} is completed.  
\subsection{Proof of Theorem \ref{Theorem3}}

a) Fixed a smooth closed $(1,1)$ form $\Omega$ such that $f^*(\theta )=\Omega +dd^cu$. If we choose another quasi-potential $u'$ so that $f^*(\theta )=\Omega +dd^cu'$ then $u-u'=\phi$, where $\phi$ is a smooth function and $dd^c\phi =0$. Condition (i) of Theorem \ref{Theorem2} is satisfied, hence the difference when we define the Monge-Ampere operator using either $u$ or $u'$ is $0$.

b) Now assume that $\{\theta\}.\{C\}=0$ for every curve in $I(f^{-1})$. Assume that we have $f^*(\theta )=\Omega +dd^cu=\Omega '+dd^cu'$, then $u-u'=\phi$ is a smooth function. Condition (ii) if Theorem \ref{Theorem2} is satisfied, therefore again we have that the difference when defining Monge-Ampere using either $\Omega$ or $\Omega '$ is $0$. 
 
\section{An example}
\label{Section1}

\subsection{The map $J_X$}
\label{Subsection1}

We consider a special Cremona transform which we will discuss in this Subsection. The material is taken from \cite{truong0}.

Let $X$ be the blowup of $\mathbb{P}^3$ along $4$ points $e_0=[1:0:0:0],e_1=[0:1:0:0],e_2=[0:0:1:0], e_3=[0:0:0:1]$;
$J:\mathbb{P}^3\rightarrow \mathbb{P}^3$ is the Cremona map $J[x_0:x_1:x_2:x_3]=[1/x_0:1/x_1:1/x_2:1/x_3]$, and let $J_X$ be the lifting of $J$ to $X$.
For $0\leq i\not= j\leq 3$, $\Sigma _{i,j}$ is the line in $\mathbb{P}^3$ consisting of points $[x_0:x_1:x_2:x_3]$ where $x_i=x_j=0$, and
$\widetilde{\Sigma _{i,j}}$ is the strict transform of $\Sigma _{i,j}$ in $X$.

Let $E_0,E_1,E_2,E_3$ be the corresponding exceptional divisors of the blowup $X\rightarrow \mathbb{P}^3$, and let $L_0,L_1,L_2,L_3$ be any lines in
$E_0,E_1,E_2,E_3$ correspondingly. Let $H$ be a generic hyperplane in $\mathbb{P}^3$, and let $H^2$ be a generic line in $\mathbb{P}^3$. Then
$H,E_0,E_1,E_2,E_3$ are a basis for $H^{1,1}(X)$, and $H^2,L_0,L_1,L_2,L_3$ are a basis for $H^{2,2}(X)$. Intersection products in complementary
dimensions are:
\begin{eqnarray*}
&&H.H^2=1,~H.L_0=0,~H.L_1=0,~H.L_2=0,~H.L_3=0,\\
&&E_0.H^2=0,~E_0.L_0=-1,~E_0.L_1=0,~E_0.L_2=0,~E_0.L_3=0,\\
&&E_1.H^2=0,~E_1.L_0=0,~E_1.L_1=-1,~E_1.L_2=0,~E_1.L_3=0,\\
&&E_2.H^2=0,~E_2.L_0=0,~E_2.L_1=0,~E_2.L_2=-1,~E_1.L_3=0,\\
&&E_3.H^2=0,~E_3.L_0=0,~E_3.L_1=0,~E_3.L_2=0,~E_3.L_3=-1.
\end{eqnarray*}
The map $J_X^*:H^{1,1}(X)\rightarrow H^{1,1}(X)$ is not hard to compute:
\begin{eqnarray*}
J_X^*(H)&=&3H-2E_0-2E_1-2E_2-2E_3,\\
J_X^*(E_0)&=&H-E_1-E_2-E_3,\\
J_X^*(E_1)&=&H-E_0-E_2-E_3,\\
J_X^*(E_2)&=&H-E_0-E_1-E_3,\\
J_X^*(E_3)&=&H-E_0-E_1-E_2.
\end{eqnarray*}
If $x\in H^{1,1}(X)$ and $y\in H^{2,2}(X)$, since $J_X^2=$the identity map on $X$, we have the duality $(J_X^*y).x=y.(J_X^*x)$. Thus from the above data,
we can write down the map $J_X^*:H^{2,2}(X)\rightarrow H^{2,2}(X)$:
\begin{eqnarray*}
J_X^*(H^2)&=&3H^2-L_0-L_1-L_2-L_3,\\
J_X^*(L_0)&=&2H^2-L_1-L_2-L_3,\\
J_X^*(L_1)&=&2H^2-L_0-L_2-L_3,\\
J_X^*(L_2)&=&2H^2-L_0-L_1-L_3,\\
J_X^*(L_3)&=&2H^2-L_0-L_1-L_2.
\end{eqnarray*}

\subsection{Proof of Theorem \ref{Theorem4}}

1) Let $\pi :Z\rightarrow X$ be the blowup at the curves $\widetilde{\Sigma _{i,j}}$. Then the induced map $J_Z:Z\rightarrow Z$ is an automorphism. If we let $g=\pi \circ J_Z:Z\rightarrow X$ then the data $(Z,\pi ,g)$ is one that we can choose to start the proof of Theorem \ref{Theorem3}.

Considering the proof of Theorem \ref{Theorem3} carefully for this choice of $Z$, it is not difficult to conclude that for a smooth closed $(1,1)$ form $\theta$ on $X$, the Monge-Ampere operator $MA(J_X^*(\theta ),\delta )$ is independent of the choice of $\delta$ if and only if the following condition is satisfied:

Condition (C): In cohomology, $\{\theta\}.\{\widetilde{\Sigma _{i,j}}\}=0$ for all $i,j$.

With the data given in Subsection \ref{Subsection1}, we can check easily that Condition (C) is satisfied if and only if $\{\theta\}=0$ is a multiple of $2H-E_0-E_1-E_2-E_3$. The latter cohomology class can be checked to be nef, by considering intersection with curves in $X$.

2) Let $X_1$ be the blowup of $\mathbb{P}^3$ at $[1:0:0:0]$. Then it is known that there is a fibration $X_1\rightarrow \mathbb{P}^2$. We then let $\Pi$ be the composition of this map with the projection from $X$ to $X_1$.

Let $h$ be a generic line in $\mathbb{P}^2$. Then in cohomology $\Pi ^*(h)=H-E_0$ where $H$ is a generic hyperplane in $\mathbb{P}^3$ and $E_0$ is the fiber over $[1:0:0:0]$. Hence 
\begin{eqnarray*}
J_X^*\Pi ^*(h)=J_X^*(H-E_0)=2H-2E_0-E_1-E_2-E_3.
\end{eqnarray*}
From the intersection on $X$, we have
\begin{eqnarray*}
J_X^*\Pi ^*(h).J_X^*\Pi ^*(h).J_X^*\Pi ^*(h)=8H^3-8E_0^3-E_1^3-E_2^3-E_3^3=-3.
\end{eqnarray*}

We then choose $\theta =\Pi ^*(\omega _{\mathbb{P}^2})$ where $\omega _{\mathbb{P}^2}$ is the Fubini-Study form on $\mathbb{P}^2$. Then $\theta \wedge \theta \wedge \theta =0$. 

On $X-I(f)$ then $MA(f^*(\theta ),\delta )=f^*(\theta \wedge \theta \wedge \theta )=0$. Therefore the support of the measure $MA(f^*(\theta ),\delta )$ is contained in $I(f)$. The total mass of this measure is the intersection in cohomology $\{f^*(\theta )\}.\{f^*(\theta )\}.\{f^*(\theta )\}=-3$.

\end{document}